\documentclass[11pt]{article}
\usepackage{amsmath,amssymb,amsfonts}
\newtheorem{thm}{Theorem}[section]
\newtheorem{cor}[thm]{Corollary}
\newtheorem{lem}[thm]{Lemma}

\newtheorem{dfn}[thm]{Definition}
\newtheorem{pro}[thm]{Proposition}

\newtheorem{ass}[thm]{Assumption}

=22cm\headheight=0cm\topskip=0cm\topmargin=0cm
\def\a{\bar {a}}\def\b{\bar {b}}
\def\x{\bar {x}}\def\y{\bar {y}}

\def\Ac{\mathcal {A}}\def\Bc{\mathcal {B}}

\def\Nn{\mathbb N}
\def\Rn{\mathbb R}

\def\M{\mathsf M}\def\N{\mathsf N}

\def\proof{{\bf Proof\ }}

\begin{document}

\begin{center}
{\LARGE{\bf Continuous integration logic}}
\vspace{5mm}

{\large{\bf Seyed-Mohammad Bagheri}}\\ \vspace{1mm}
{\footnotesize Department of mathematics, Tarbiat-Modares University,\\
Tehran, Iran, P.O. Box 14115-134;\\
e-mail: bagheri@modares.ac.ir}
\vspace{1mm}

{\large{\bf Massoud Pourmahdian}}\\ \vspace{1mm}
{\footnotesize School of Mathematics, Amirkabir University of
Technology, \\ Hafez avenue 15194, Tehran, Iran;\\
e-mail: \ pourmahd@aut.ac.ir}
\vspace{5mm}

\end{center}

\bigskip


\begin{abstract} We combine continuous and integral logics and found
a logical framework for metric measure spaces equipped with a
family of continuous relations and operations.
We prove the ultraproduct theorem and deduce compactness and other
usual results. We also give applications of the compactness theorem
in metric measure theory.
\end{abstract}

\section{Introduction}
Classical model theory is usually described as the study of algebraic structures by
logical methods. The efficiency of the program has always been a good reason
for introducing other model theories. One of the first such model theories is
Chang and Keisler's continuous model theory \cite{Chang-Keisler2}.
Several other variants was then introduced for special purposes
including Banach space model theory \cite{Henson-Iovino}, model theory of
probability structures (or probability logic) \cite{Hoover, Keisler1, Fajardo, integration}
and model theory of metric structures (also called continuous logic) \cite{BBHU}.
In all these logics ultraproduct construction plays the central role.
It is used to prove the compactness theorem which is the cornerstone of model theory.
Apart from logical considerations, the ultraproduct construction is a powerful and
flexible method for constructing new structures. Ultraproducts of topological spaces,
metric spaces, measure spaces and many other kinds of structures has been defined
and studied. The ultraproduct of a family of metric measure spaces
is defined analogously and is a metric measure space again.
This fact can be used to introduce a relevant model theory
for structures equipped with a metric and a measure.

The purpose of this paper is to combine continuous logic with
probability logic and introducing a model theory for metric structures
equipped with a compatible measure (i.e. a topological measure).
Here, by probability logic is meant integration logic, i.e. one which
uses integration as a quantifier. Topological measures on metric spaces
are usually called metric measure spaces in the literature.
Most arguments in the formation of this logic is similar to continuous
logic so that the existing proofs work in this setting as well.
There are however aspects which are special to the present context
and relate measure to the metric. The most important of these aspects
is the ultraproduct construction where the measure is defined in such
a way that the result is a topological measure.
There are also aspects which are completely new.
We mostly presents the essential features of the combined logic
which will be called continuous integration (or integral) logic in the paper.
We also give some applications of the compactness theorem in metric measure theory.
In particular, we give a short proof of the Riesz representation theorem
stating that every positive linear functional on $C(X)$,
where $X$ is a compact metric space, is an integral.

The paper is organized as follows. In the next section we gather up the
essential background needed for the foundation of continuous integral logic.
In section \ref{integral-metric} we review basics of continuous logic and
integral logic. For the sake of comparison, we sometimes refer to the
first as metric logic and to the second as measure logic.
In section \ref{metric-integral} we discuss the combination of these logics.
In the last section we will give two applications of the compactness theorem.

\section{Preliminaries from measure theory}
In this section we review some basic facts from measure theory.
A measure on a Boolean algebra $\Bc$ of subsets of $M$ is a map
$\mu:\Bc\rightarrow[0,\infty]$ such that $\mu(\emptyset)=0$ and for any
countable sequence $A_k\in \Bc$ of disjoint sets for which $\cup_kA_k\in\Bc$,
$$\mu(\cup_kA_k)=\sum_k\mu(A_k).$$
If $\Bc$ is a $\sigma$-algebra, $\mu$ is called a measure.
An outer measure on $M$ is a map $\mu^*:\mathcal P (M)\rightarrow[0,\infty]$
such that (i) $\mu^*(\emptyset)=0$
(ii) $\mu^*(A)\leqslant\mu^*(B)$ whenever $A\subseteq B$
(iii) $\mu^*(\cup_kA_k)\leqslant\sum_k\mu^*(A_k)$ for any sequence $A_k$.
If $\mu^*$ is an outer measure on $X$, a set $E\subseteq X$ is called
$\mu^*$-measurable if for every $A$ $$\mu^*(A)=\mu^*(A\cap E)+\mu^*(A-E).$$

\begin{thm} \emph{(Carath\'eodory extension theorem)}
Let $\mu$ be a measure on a Boolean algebra $\Ac$ of subsets of
$M$. Then $\mu$ has an extension $\bar\mu$ to the $\sigma$-algebra
of $\mu^*$-measurable sets. If $\mu$ is $\sigma$-finite,
$\bar\mu\upharpoonright\sigma(\Ac)$ is unique.
\end{thm}

Let $M$ be a metric space. The distance between two sets $X,Y\subseteq M$
is $$d(X,Y)=\inf \{d(x,y):\ x\in X, y\in Y\}.$$
$X,Y$ are called separated if $d(X,Y)>0$.
An outer measure $\mu^*$ on a metric space $(M,d)$ is said to be a
{\em metric outer measure} if for every separated sets $X,Y$ one has that
$$\mu^*(X\cup Y)=\mu^*(X)+\mu^*(Y).$$

\begin{pro} \emph{(\cite{Bruckner} Th. 3.7)} \label{Bruckner}
Let $\mu$ be an outer measure on a metric space $M$. Then every Borel subset
of $M$ is $\mu^*$-measurable if and only if $\mu$ is a metric outer measure.
\end{pro}

A measure $\mu$ on a topological space $M$ is called a {\em Borel measure} if
all Borel sets are measurable.
A measure $\mu$ on a metric space $M$ is called a {\em metric measure} if for each
metrically separated $X,Y$, one has that $\mu(X\cup Y)=\mu(X)+\mu(Y)$.
It is a fact that a measure $\mu$ on a metric space $X$ is a metric measure
if and only if it is a Borel measure.

\begin{pro} \label{Borel extension}
Let $(M,d)$ be a metric space and $\mu$ a finite measure on a
$\sigma$-subalgebra $\Sigma$ of the $\sigma$-algebra of Borel sets of $M$.
Then $\mu$ has an extension to a metric measure $\hat\mu$ on $M$.
\end{pro}
\proof For $\delta>0$ and set $E\subseteq M$, set
$$\nu^*_\delta(E)=\inf \big\{\sum_k\mu(C_k):\ E\subseteq\cup_kC_k,\ diam(C_k)\leqslant\delta\big\}$$
$$\nu^*(E)=\lim_{\delta\rightarrow0}\nu^*_\delta(E).$$
Then one can readily see that $\nu^*$ defines a metric outer measure.
Therefore, $\nu^*$ induces a measure $\nu$ on the family of Borel sets.
Moreover, as $\mu$ is a measure, for any set $E\in\Sigma$, $\mu(E)=\nu(E)$. $\square$
\bigskip

A Borel measure $\mu$ on a topological space $M$ is said to be \emph{regular} if for
every Borel $A\subseteq M$,
$$\mu(A)=\sup\{\mu(C):\ C\subseteq A, \ C\ \text{closed}\}=
\inf\{\mu(U):\ A\subseteq U, \ U\ \text{open}\}.$$
Any Borel measure on a metric space is regular. If it is complete and separable, it is Radon,
i.e. measures of subsets are approximated from below by compact sets
(\cite{BogachevII} \S 7).
A Borel measure $\mu$ on $M$ is called \emph{$\tau$-additive} if for any upward directed
family $\{U_i\}_{i\in I}$ of open sets $$\mu(\cup_{i}U_i)=\sup_{i} \mu(U_i).$$
Every Borel measure on a separable metric space is a $\tau$-additive (see \cite{BogachevII}).
The \emph{weight} of a metric space $(M, d)$ is the minimal cardinality
of a topology base in $M$ (and also the minimal cardinality $\kappa$ with the
property that every set $S\subseteq M$ with $\inf_{x,y\in S, x\neq y} d(x, y)>0$
is of cardinality at most $\kappa$.
The support of a topological measure $\mu$ on $M$ is the smallest closed
set $E$ such that $\mu(M-E)=0$. Every $\tau$-additive measure space has a support.

\begin{pro} \label{tau-additive} \emph{(\cite{BogachevII} 7.2.10.)}
The weight of a metric space $(M,d)$ is a nonmeasurable cardinal precisely when
every Borel measure on $M$ is $\tau$-additive (and then is Radon if X is complete).
An equivalent condition: every Borel measure on $M$ has support.
\end{pro}

\begin{cor} \label{uniqueness} \emph{(\cite{BogachevII} 7.2.3.)}
Let two $\tau$-additive measures $\mu$ and $\tau$ on a space $M$ coincide on all sets
from some class $U$ that contains a base of the topology in $M$ and is closed with
respect to finite intersections. Then $\mu=\tau$.
\end{cor}

Given $\tau$-additive Borel probability spaces $M$, $N$, there is a
unique $\tau$-additive Borel measure on $M\times N$ extending the
product measure. It is called the \emph{$\tau$-additive product}.

\begin{pro} \label{TauFubini} \emph{(\cite{Fremlin} 417)}
Let $(M,\mu)$, $(N,\nu)$ be $\tau$-additive Borel probability spaces.
If $B\subseteq M\times N$ is Borel then for each $a\in M$, $B_a$ is Borel
and the function $\nu(B_x)$ is Borel measurable. If $f:M\times N\rightarrow \Rn$
is bounded and continuous then $x\mapsto\int f(x,y)d\mu$ is continuous and
$$\iint fd\mu d\nu=\iint fd\nu d\mu= \int f d(\mu\times\nu).$$
\end{pro}

Let $M$ be a topological space, $\mathcal D$ an ultrafilter
on an index set $I$ and $\{x_i\}_{i\in I}$ a family of points in $M$.
Then one sets $\lim_{\mathcal D}x_i=x$ if for each open $U$ containing
$x$ one has that $\{i: x_i\in U\}\in\mathcal D$.
It is a basic fact that $X$ is compact Hausdorff if and only if every such
sequence has a $\mathcal D$-limit.

Let $(M,\Ac,\mu)$ be a measure space and $\pi:M\rightarrow N$ a map.
The \emph{push forward measure} $\pi^*\mu$ on $N$ is defined by
$$\ \ \ \ \pi^*\mu(B)=\mu(\pi^{-1}(B)) \ \ \ \ \ \ \ \ \ \ B\subseteq N\
\mbox{and}\ \pi^{-1}(B)\in\Ac.$$
By the change of variable formula, if $f:N\rightarrow\Rn$ is $\pi^*\mu$-measurable
then $f\circ\pi$ is $\mu$-measurable and $\int_M f\circ\pi\ d\mu=\int_N f\ d(\pi^*\mu)$.

Let $(N,\Bc,\mu)$ be a measure space. The outer measure $\mu^*$ on $N$ is defined by
$$\mu^*(X)=\inf\{\mu(A)|\ X\subseteq A\in\Bc\}.$$ If $M\subseteq N$ then
$\Bc_M=\{A\cap M|\ A\in \Bc\}$ is a $\sigma$-algebra
and the restriction of $\mu^*$ to $\Bc_M$, denoted by $\mu_M$, is a
measure on $M$. In fact, elements of $\Bc_M$ are
$\mu^*|_M$-measurable. $\mu_M$ is called the {\em subspace measure}
on $M$. A measurable envelope for $M$ is a measurable set $E\in\Bc$
such that $M\subseteq E$ and $\mu(E\cap A)=\mu^*(M\cap A)$ for any
$A\in\Bc$. Every $M\subseteq N$ of finite outer measure has an
envelope. In fact, every $E\in\Bc$ containing $M$ with
$\mu(E)=\mu^*(M)$ is a measurable envelope for $M$ (see
\cite{Fremlin}, 132E). If $f:N\rightarrow \Rn$ is measurable, by
$\int_M f$ is meant $\int_M (f|_M)d\mu_M$.

\begin{pro} \label{subspacemeasure} \emph{(\cite{Fremlin}, 214)}
Let $(N,\Bc,\mu)$ be a measure space, $M\subseteq N$ and $f$ an
integrable function defined on $N$.

(i) $f|_M$ is $\mu_M$-integrable and $\int f\geqslant\int_M f$ if\ $f\geqslant0$.

(ii) If either $M$ is of full outer measure in $N$ or $f$ is zero
almost everywhere on $N-M$, then $\int_Mf=\int_Nf$.

(iii) If $E\in\Bc$ is a measurable envelope of $M$ then $\int_Mf$ is
equal to $\int_Ef=\int_Nf\cdot {\mathcal X}_E$.
\end{pro}

A modulus of uniform continuity is a function $\Delta:\Rn^+\rightarrow\Rn^+$.
A map $f:(M,d)\rightarrow (M',d')$ is said to be uniformly continuous with
modulus $\Delta$ if $d(x,y)<\Delta(\epsilon)$ implies $d'(f(x),f(y))\leqslant\epsilon$.

\section{Integration logic and continuous logic} \label{integral-metric}
A {\em basic language} is a usual first order language consisting of constant,
function and relation symbols. To each relation symbol $R$ (resp function symbol $F$)
is assigned a natural number $n_R\geqslant 1$ (resp $n_F\geqslant 1$) called its arity.
Also, to each relation symbol $R$ is assigned a real number $\flat_R\geqslant0$ called
its (uniform) bound.
In the next sections, we will put further conditions on the languages to obtain
metric, measure or metric-measure languages.
The set of real numbers is always is used as value space of logic.
Logical symbols consist of the connectives and quantifiers.
The primitive connectives used in this paper are $+,\wedge$ and
scalar product $r\cdot$ for each $r\in\Rn$. Other connectives such as
$-, \vee$ and absolute value $|\cdot |$ are obtained by combining them in the
obvious way. The needed quantifiers depend on the logic.
In integration logic it is $\int$, in metric logic is `$\sup$' and in metric-integration
logic are both. A {\em basic $L$-structure} is a nonempty set $M$ equipped with

- for each constant symbol $c$, an element $c^{M}$

- for each relation symbol $R$, a function $R^{M}:M^{n_R}\rightarrow[-\flat_R,\flat_R]$

- for each function symbol $F$ a function $F^{M}:M^{n_F}\rightarrow M$.
\bigskip

$L$-terms are defined in the usual way, i.e. constant symbols and variables are terms and
if $F$ is a $n$-ary function symbol and $t_1,...,t_n$ are terms then $F(t_1,...,t_n)$ is a term.
For each basic $L$-structure $M$ and term $t(x_1,...,x_n)$ the function
$t^M:M^n\rightarrow M$ is defined in the obvious way.

\subsection{Integration logic} \label{integration logic}
In this section we review some basic facts from Integration logic.
For more details see \cite{integration}, \cite{Hoover} and \cite{Keisler1}.
An {\em measure language} is a basic language containing a distinguished binary relation
symbol $e$ for equality (with $\flat_e=1$) and equipped with a second order symbol $\mu$ for measure.
The connectives are as stated above and the only quantifier symbol is integration $\int$.
Let $L$ be a measure language. $L$-terms are defined in the usual way.
Formulas and their bounds are defined inductively as follows:

- $1$ is a formula with bound $1$

- if $R$ is a $n$-ary relation symbol and $t_1,...,t_n$ are terms, $R(t_1,...,t_n)$
is an atomic formula with bound $\flat_R$

- if $\phi,\psi$ are formulas and $r\in\Rn$, then $r\phi$ is a formula with
bound $|r|\flat_\phi$ and $\phi+\psi$, $\phi\wedge\psi$ are formulas with bound
$\flat_\phi+\flat_\psi$

- if $\phi$ is a formula and $x$ is a variable, then $\int\phi dx$ is a formula
with bound $\flat_\phi$.
\bigskip

The notion of free variable is defined in the obvious way.
Every formula can be displayed in the form $\phi(\x)$ where $\x$ is the list of its
free variables. A sentence is a formula without free variables.
Expressions of the form $\phi=\psi$ or $\phi\leqslant\psi$ are called \emph{statements}.
If $\phi,\psi$ are sentences, the corresponding statements are called \emph{closed statements}.

\begin{dfn} \label{graded}
\em{A {\em (graded) measure $L$-structure} is a basic $L$-structure $M$
equipped for each $n$ with a measure $(\Bc_n,\mu_n)$ on $M^n$
 such that the following conditions hold:

\begin{enumerate}
\item{$\mu_1(M)\leqslant 1$ and for all $m,n$,
$\mu_{m+n}$ is an extension of the product measure $\mu_m\times\mu_n$.}

\item{Each $\mu_n$ is invariant under the permutations of variables.}

\item{For every terms $t_1(\x),...,t_k(\x)$, the map
$\x\mapsto (t_1^M(\x),...,t_k^M(\x))$ is measurable.
Every $R^M(x_1,...,x_n)$ is $\mu_n$-measurable.}

\item{The Fubini property holds: if $B$ is $\mu_{m+n}$-measurable then}

\begin{itemize}
\item{for all $\a\in M^m$, $B_{\a}=\{\b\in M^n|\ (\a,\b)\in B\}$
is $\mu_n$-measurable,}

\item{the function $\x\mapsto\mu_n(B_{\x})$ is
$\mu_m$-measurable,}

\item{$\int \mu_n(B_{\x})d{\mu_m}=\mu_{m+n}(B)$.}
\end{itemize}
\end{enumerate}}
\end{dfn}

Note that the diagonals are usually non-measurable in the product measures
so that $\mu_n$ is generally a proper extension of the product measures.
We will denote measure structures by $\M,\N$ etc. Let $\M$ be a $L$-structure.
Formulas are interpreted inductively as follows:
\bigskip

- $1^{\sf M}=1$

- $(R(t_1(\x),...,t_k(\x)))^{\sf M}(\a)=R^{\sf M}(t^{\sf M}_1(\a),...,t^{\sf M}_k(\a))$

- $(r\phi+s\psi)^{\sf M}=r\phi^{\sf M}+s\psi^{\sf M}$

- $(\phi\wedge\psi)^{\sf M}=\phi^{\sf M}\wedge\psi^{\sf M}$

- $(\int\phi(x)dx)^{\sf M}=\int\phi^{\sf M}(x)d\mu$
\bigskip

The following lemma is easily proved by induction on the complexity of formulas.

\begin{lem}
For every formula $\phi(\x)$, $\phi^{\sf M}$ is a real valued
measurable function bounded by $\flat_\phi$.
\end{lem}

The Fubini property implies that Fubini's theorem holds for every formula:
$$\int\int\phi^{\sf M}dxdy=\int\int\phi^{\sf M}dydx.$$

One other reason for considering measure structures in the graded form
is that the ultraproduct construction works for them.
Let $\M_i=(M_i,\Bc_i^n,\mu_{in})$, $i\in I$, be an indexed family of
measure $L$-structures and $\mathcal D$ an ultrafilter over $I$.
Let $M=\prod_\mathcal DM_i$ be the set theoretic ultraproduct of the family.
So, elements of $M$ are equivalence classes of the relation defined on
$\prod_{i\in I}M_i$ by setting $(x_i)\sim(y_i)$ if $\{i|\ x_i=y_i\}\in\mathcal D$.
The equivalence class of $(a_i)$ is denoted by $[a_i]$.
First, we put a basic $L$-structure on $M$ by interpreting the symbols of $L$ as follows:

- $c^M=[c^{M_i}]$

- $F^M([a_i^1],...,[a_i^n])=[F^{M_i}(a_i^1,...,a_i^n)]$

- $R^M([a_i^1],...,[a_i^n])=\lim_{\mathcal D} R^{M_i}(a_i^1,...,a_i^n)$.

We also define a measure on any $M^n$. First assume $n=1$.
If $A_i\subseteq M_i$ is $\mu_{i1}$-measurable, the set
$$[A_i]=\big\{[a_i]\ : \ \{i: a_i\in A_i\}\in\mathcal D \ \big\}\subseteq M$$
is called an ultrabox in $M$. Ultraboxes in $M$ form a Boolean algebra.
Moreover, $[A_i]\subseteq[B_i]$ if and only if $\{i|\ A_i\subseteq B_i\}\in\mathcal D$.
Define a real valued map on ultraboxes by setting
$$\mu_1([A_i])=\lim_{\mathcal D}\mu_{i1}(A_i).$$
It is not hard to see that

\begin{lem}
$\mu_1$ is a measure on the Boolean algebra of ultraboxes in $M$.
\end{lem}

Identifying $M^n$ with $\prod_{\mathcal D}M_i^n$ in the natural way, one can
similarly define a measure $\mu_n$ on the Boolean algebra of ultraboxes in $M^n$.
By the Carath\'eodory extension theorem, each $\mu_n$ extends to a
unique measure on the $\sigma$-algebra $\Ac^n$ generated by the
ultraboxes of $M^n$. We continue denoting this measure
by $\mu_n$.
The interpretations of symbols of the language are measurable with respect to these measures.
For example, assume $R$ is a unary relation symbol.
Let $J=(r,s)$ be an open interval and $J=\cup_n I_n$ where $\{I_n\}_n$
is an increasing sequence of closed intervals. One checks easily that
$$(R^\M)^{-1}(J)=\bigcup_n[(R^{\M_i})^{-1}(I_n)].$$
Indeed, it can be shown that

\begin{pro}\label{u1}
$\M=(M,\Ac^n,\mu_n)$ is a measure $L$-structure.
\end{pro}

The following lemmas are crucial for the proof of integral {\L}o\'{s} theorem.

\begin{lem} \label{exchange1}
For each $i\in I$, let $a_i^k$ be a sequence of real numbers tending to $a_i$
where $|a_i|\leqslant r$. Assume these sequences are uniformly convergent,
i.e. for each $\epsilon>0$ there is a $n_\epsilon$ such that for any $i\in I$ and
$k\geqslant n_\epsilon$, $|a_i^k-a_i|<\epsilon$ .
Then $\lim_k \lim_{\mathcal D}a_i^k=\lim_{\mathcal D}a_i$.
\end{lem}

\begin{lem} \label{exchange2}
For each $i\in I$, let $f_i:M_i\rightarrow\Rn$ be a bounded by
$\beta\geqslant0$ measurable function. Then
$\int_{M}\lim_{\mathcal D} f_i=\lim_{\mathcal D}\int_{M_i}f_i$.
\end{lem}

\begin{thm} \label{IntegralLos} \emph{(Ultraproduct theorem)}
For any formula $\phi(\x)$ and $[a_i^1],...,[a_i^n]$
$$\phi^\M([a_i^1],...,[a_i^n])=
\lim_{\mathcal D}\phi^{\M_i}(a_i^1,...,a_i^n).$$
\end{thm}
\proof The claim is proved by induction on the complexity of $\phi$.
The atomic and connective cases are obvious.
Consider the case $\psi=\int\phi(\x,y)dy$.
Then by Lemma \ref{exchange2} and induction hypothesis
$$\int\phi^\M(\a,y)dy=\int \lim_{\mathcal D}\phi^{\M_i}(\a_i,y)dy=
\lim_{\mathcal D}\int\phi^{\M_i}(\a_i,y)dy=\lim_{\mathcal D}\psi^{\M_i}(\a_i).\ \square$$

\subsection{Continuous (metric) logic}
A {\em metric language} is basic language $L$ containing a distinguished binary symbol
$\rho$ for metric and equipped for each relation symbol $R$ (resp function symbol $F$)
with a modulus of uniform continuity $\Delta_R$ (resp $\Delta_R$).
We always assume that $\flat_\rho=1$ and $\Delta_\rho=id$.
Logical symbols consist of the connectives $+,\wedge, r\cdot$
as before and the quantifier $\sup$. We also set $\inf=-\sup-$.

Let $L$ be a metric language.
The collection of $L$-terms and their modulus of uniform continuities are defined inductively.
In particular, the modulus of continuity of $F(t_1,...,t_n)$
is $\min_k \Delta_{t_k}(\Delta_F(\epsilon))$.
The collection of $L$-formulas with their uniform bounds and modulus of uniform
continuities are defined inductively as follows:

- $1$ is an atomic formula with uniform bound $1$ and
modulus of uniform continuity $0$

- If $R$ is a $n$-ary relation symbol and $t_1,...,t_n$ are terms,
then $R(t_1,...,t_n)$ is an atomic formula with uniform bound $\flat_R$ and
modulus of uniform continuity $\min_k \Delta_{t_k}(\Delta_R(\epsilon))$

- If $\phi$ is a formula and $r\neq0$, then $r\phi$ is a formula with bound
$|r|\phi$ and modulus of uniform continuity $\Delta_\phi({\epsilon\over |r|})$;
$0\phi$ is a formula with bound $0$ and modulus of uniform continuity $0$

- If $\phi$ and $\psi$ are formulas, then so are $\phi+\psi, \ \phi\wedge\psi$
with bound $\flat_\phi+\flat_\psi$ and modulus of uniform continuity
$\min \{\Delta_\phi({\epsilon\over 2}), \Delta_\psi({\epsilon\over 2})\}$

- If $\phi$ is a formula then so is $\sup_x\phi$ with bound $\flat_\phi$
and modulus of uniform continuity $\Delta_\phi$.
\bigskip

Let $L$ be a metric language. A \emph{metric $L$-structure} $\M$ is a basic $L$-structure
$M$ equipped with a metric $\rho^{\M}$ of diameter at most $1$ such that
every $R^{\M}$ (resp $F^{\M}$) is uniformly continuous with modulus $\Delta_R$
(resp $\Delta_F$) where we put the maximum metric on the Cartesian powers.
Let $\M$ be a metric $L$-structure, $\phi(\x)$ be a formula and $\a\in M$.
Then $\phi^{\M}(\a)$ is defined similar to the integration logic. In particular,
$(\sup_x\phi)^{\M}=\sup_{a\in M}\phi^{\M}(a)$.
The following proposition is easily proved by induction on the complexity of formulas.

\begin{pro}
$\phi^{\M}(\x)$ is uniformly continuous with modulus $\Delta_\phi$.
Moreover, $|\phi^{\M}|\leqslant\flat_\phi$.
\end{pro}

Now we define the ultraproduct construction for metric structures.
Let $L$ be a metric language, $\M_i=(M_i,\rho_i)$, $i\in I$, be an
indexed family of metric $L$-structures and $\mathcal D$ be an ultrafilter over $I$.
Let $M=\prod_\mathcal DM_i$ be the set theoretic ultraproduct of the family
and put a basic $L$-structure on $M$ as in the previous subsection.
In particular, for $a=[a_i], b=[b_i]$ one sets
$$\rho(a,b)=\lim_{\mathcal D} \rho_i(a_i,b_i).$$
Then $\rho$ is a pseudometric on $M$ and so $\rho(a,b)=0$ defines an equivalence relation on $M$.
We denote the class of $a$ by $\hat a$. Then $\rho$ induces a metric on
the quotient set $\hat M$ which we denote by $\hat\rho$.
Note that the uniform continuity of the relations $R^{M_i}$ with resect to
the modulus $\Delta_R$ implies that $R^{M}$ induces a well-defined function on $\hat M$.
Similarly, $F^{M}$ induces a well-defined function on $\hat M$.
We denote the resulting metric $L$-structure by $\hat\M$.

\begin{thm} \label{MetricLos} \emph{(Ultraproduct theorem)}
For any formula $\phi(\x)$ and $a=[a_i^1], b=[b_i^n]$
$$\phi^{\hat\M}(\hat a,\hat b,...)=
\lim_{\mathcal D}\phi^{\M_i}(a_i,b_i,...).$$
\end{thm}

\section{Continuous integration logic} \label{metric-integral}
In this section we combine integration and metric logics.
We fix a set theoretic assumption which facilitates technical details.

\begin{ass} \label{assumption}
Every Borel measure on a metric space is $\tau$-additive.
\end{ass}

By Proposition \ref{tau-additive}, if there is no measurable cardinal in the universe
(this is in particular true if $V=L$ holds) then this assumption holds.
We recall also that every measurable cardinal is inaccessible and it is well-known that
the consistency of ZFC+$\exists$(inaccessible cardinal) is not provable.
Indeed, this assumption is just for convenience and
all what follows can be done with a lot of further complication.

\subsection{Syntax and semantics}
A {\em metric-measure language} is a metric language equipped with a measure symbol $\mu$.
Logical symbols consist of the connectives $+,\wedge, r\cdot$ as before, and the
quantifiers $\sup$ and $\int$.
Formulas are defined as in metric logic with a further formula making rule:

- if $\phi$ is a formula with bound $\flat_\phi$ and modulus of continuity
$\Delta_\phi$, and $x$ is a variable, then $\int\phi dx$ is a
formula with bound $\flat_\phi$ and modulus of continuity $\Delta_\phi$.

\begin{dfn}
{\em A {\em metric-measure $L$-structure} is a metric structure $(M,\rho)$ in $L$
equipped with a Borel measure $\mu$ on $M$ such that $\mu(M)\leqslant1$.}
\end{dfn}

We put the maximum metric and also the $\tau$-additive product measure on every $M^n$.
Using Proposition \ref{TauFubini} we can easily check that the conditions of Definition
\ref{graded} hold for metric-measure structures so that

\begin{pro}
Every metric-measure $L$-structure is graded.
\end{pro}

Let $L$ be a metric-measure language and $\M$ a structure in $L$.
Let $\phi(\x)$ be an $L$-formula and $\a\in M$. Then $\phi^{\M}(\a)$ is defined
by induction on the complexity of $\phi$ as in metric and integration logics.
In particular, if $\phi^{\M}(\x,y)$ is defined, then

- $(\sup_y\phi)^{\M}(\a)=\sup_{b\in M}\phi^{\M}(\a,b)$

- $(\int\phi dy)^{\sf M}(\a)=\int\phi^{\sf M}(\a,y)dy$.
\bigskip

The following proposition is easily proved by induction on the complexity of formulas.

\begin{pro}
$\phi^{\M}(\x)$ is uniformly continuous with modulus $\Delta_\phi$
and $|\phi^{\M}|\leqslant\flat_\phi$. In particular, it is measurable.
\end{pro}

We now describe the ultraproduct construction in the framework of
metric-integration logic.
Let $L$ be a metric-measure language, $(M_i,\mu_i,\rho_i,...)_{i\in I}$ be
an indexed family of $L$-structures and $\mathcal D$ be an ultrafilter on $I$.
Let $M=\prod_{\mathcal D} M_i$ and let $(\hat M,\hat\rho)$ be the resulting metric structure
as defined in metric logic. Recall that $\hat M$ is a quotient of $M$.
We wish to put a measure on $\hat M$ turning it to a metric-measure $L$-structure.
For this purpose, we first put the ultraproduct measure $\mu$ on $M^1$ as defined in integration logic.
Let $\pi:M\rightarrow\hat M$ be the quotient map and $\pi^*\mu$ be the push forward
measure on $\hat M$.
Note that in $M$ with its pseudometric $\rho$ we have that
$$B([a_i];s):=\big\{[x_i]: \rho([a_i],[x_i])<s\big\}=\bigcup_{r<s} [B(a_i;r)].$$
This shows that every open ball in $(\hat M,\hat\rho)$ is $\pi^*\mu$-measurable.
Thus, by Propositions \ref{tau-additive} and \ref{uniqueness}, the restriction of
$\pi^*\mu$ to the $\sigma$-algebra generated by the balls extends uniquely to a Borel
measure on $\hat M$ which we denote by $\hat\mu$.
Thus, $\hat\M=(\hat M,\hat\rho,\hat\mu)$ is a metric-measure $L$-structure.

\begin{lem} \label{exchange3}
For each $i\in I$, let $f_i:M_i\rightarrow\Rn$ be bounded by $\beta\geqslant0$ and
uniformly continuous with modulus $\Delta$. Let $f$ be the function induced
by $\lim_{\mathcal D}f_i$ on $\hat M$. Then
$\int_{\hat M} f=\lim_{\mathcal D}\int_{M_i}f_i$.
\end{lem}
{\bf Proof} By Lemma \ref{exchange2} we must show that $\int_{\hat M} f=\int_{M}\lim_{\mathcal D} f_i$.
Clearly, this equality is a special case of the change of variable formula if we can show that
$f$ is $\pi^*\mu$-measurable. But, this latter a consequence of uniform continuity of all
$f_i$'s with respect to $\Delta$. $\square$
\bigskip

\begin{thm} \emph{(Ultraproduct theorem)} \label{Los theorem}
For each formula $\phi(\x)$ in the metric-measure language $L$ and $a_1=[a_i^1],...,a_n=[a_i^n]$
$$\phi^{\hat\M}(\pi a_1,...,\pi a_n)=\lim_{\mathcal D}\phi^{\M_i}(a_i^1,...,a_i^n).$$
\end{thm}
\proof The proof is done by induction on the complexity of formulas. The integration
step is by Lemma \ref{exchange3} and the supremum step is as in theorem \ref{MetricLos}.
$\square$
\bigskip

Note that if every $M_i$ is complete then so is $\hat M$.
The notions of elementary embedding, elementary equivalence etc are defined in the obvious way.
The first consequence of the ultraproduct theorem is the compactness theorem.
In some applications, it is better to use an approximate version of this theorem.
We say a finite set of closed statements $\sigma_1=0,\ldots,\sigma_k=0$ is
\emph{approximately satisfiable} if for each $\epsilon>0$ there is a model $M$
such that $|\sigma_i^M|\leqslant\epsilon$, $i=1,...,k$.
Using a suitable nonprincipal ultrafilter on $\Nn$ one can show that an
approximately satisfiable finite set of statements is satisfiable. We have then

\begin{thm} \emph{(Compactness)}
Every finitely satisfiable set of closed statements is satisfiable.
Every finitely approximately satisfiable set of closed statements is satisfiable.
\end{thm}

The following propositions are also proved as in continuous logic.

\begin{pro} \label{axiomatization} \emph{(Axiomatizability)}
A class $\mathcal K$ of $L$-structures is an elementary class if and
only if it is closed under elementary equivalence and ultraproduct.
\end{pro}

The elementary diagram of $M$, denoted $ediag(M)$, is the set
of all $L(M)$-sentences $\phi(\a)$ such that $\phi^{M}(\a)=0$.
It is clear that $N\vDash ediag(M)$ if and only if $M\preceq N$.

\begin{pro} \emph{(Elementary AP and JEP)}
Let $f:M_0\rightarrow M_1$ and $g:M_0\rightarrow M_2$ be elementary embeddings.
Then there are $N$ and elementary embeddings $f':M_1\rightarrow N$ and
$g':M_2\rightarrow N$ such that $f'f=g'g$. Similarly any two elementarily
equivalent structures are elementarily embedded in a third structure.
\end{pro}

Now we prove the union of elementary chains theorem.

\begin{pro} \emph{(Union of chains)}
Let $${\sf M}_0\preceq {\sf M}_1\preceq\cdots {\sf M}_\alpha\preceq\cdots \ \ \ \ \ \alpha\in\kappa$$
be an elementary chain of $L$-structures.
Then there is a $L$-structure ${\sf M}=\cup_\alpha {\sf M}_\alpha$ such that
${\sf M}_\alpha\preceq {\sf M}$ for any $\alpha$.
\end{pro}
\proof: $M=\bigcup M_\alpha$ is obviously a metric structure.
Now we define a metric measure on $M$ in order to obtain a metric-measure structure.
For this purpose, it is enough to define a metric outer measure on $M$.
Suppose $\mu^*_\alpha$ is the metric outer measure associated to $\mu_\alpha$.
Then for each $E\subset M$ set
$$\mu^*(E)=\sup_\alpha \mu^*_\alpha (E\cap M_\alpha).$$
An easy calculation shows that $\mu^*$ is a metric outer measure.
Hence $\mu^*$ induces a measure $\mu$ on the set of Borel subsets of $M$.
We want to show that for any $\alpha$, $M_\alpha\preceq M$.
We should proceed the proof by induction on the complexity of $L$-formulas.
The claim holds for atomic formulas obviously.
Assume it holds for $\phi(x,\y)$.
Hence for any $a,\b\in M_\alpha$, we have that $\phi^{\sf M}(a,\b)=\phi^{{\sf M}_\alpha}(a,\b)$.
Notice that for any $r<s$,
$$\mu\{x: r< \phi^{\sf M}(x,\b)\leqslant s\}=\mu_\alpha\{x: r< \phi^{{\sf M}_\alpha}(x,\b)\leqslant s\}.$$
Thus $$\int_{M} \phi^{\sf M}(x,\b) dx=\int_{M_\alpha} \phi^{{\sf M}_\alpha} (x,\b)dx.$$
The induction step for $\sup$ and the connectives are similar to metric logic.
$\square$

\subsection{Substructure}
Let $M,N$ be metric-measure $L$-structures and $M$ be a subset of $N$.
We say $M$ is a substructure of $N$ if the interpretation of every symbol on
$M$ equals the restriction to $M$ of the corresponding symbol on $N$.
Equivalently, for every quantifier-free $\phi$ one must have that
$\phi^N|_M=\phi^M$. If this equality holds for every $L$-formula we say
$M$ is an elementary substructure of $N$ and write $M\preceq N$.

\begin{pro} \label{Tarski-Vaught}\emph{(Tarski-Vaught)}
Assume $M\subseteq N$. Then $M\preceq N$ if and only if
for each $\phi(x)$ with parameters in $M$,

\emph{(i)} $\sup_{x\in M}\phi^N(x)=\sup_{x\in N}\phi^N(x)$.

\emph{(ii)} $\{x\in M|\ \phi^N(x)>0\}$ is $\mu_M$-measurable and has the same measure.
\end{pro}
\proof The `only if' part is clear. Let us prove the `if' part.
Assume the mentioned conditions hold.
We show by induction that for each $\phi(\x)$ and $\a\in M$,
$\phi^M(\a)=\phi^N(\a)$. Atomic, connective and $\sup$ cases are easy.
Let consider the integral case. Assume the claim holds for $\phi(\x,y)$ and fix $\a\in M$.
The assumption of the proposition implies that for any $r,s$, the expression
$r<\phi(\a,y)\leqslant s$ defines sets of the same measure in $M$ and $N$.
Partitioning the interval $[-\flat_\phi,\flat_\phi]$ to small enough intervals
of the form $(r,s]$, we can define a sequence $\xi_k$ of simple $\mu_{N}$-measurable
functions on $N$ tending to $\phi^\N(\a,y)$.
Let $\eta_k$ be the restriction of $\xi_k$ to $M$. Then $\eta_k$ is
measurable on $M$ and moreover $$\int_M\eta_kd\mu_{M}=\int_N\xi_kd\mu_{N}.$$
Therefore, $$\int_M\phi^M(\a,y)dy=\lim_k\int_M\eta_k(y)dy=
\lim_k\int_N\xi_k(y)dy=\int_N\phi^N(\a,y)dy. \ \square$$

It is natural to ask whether every metric-measure structure can be completed.
We show that this is done in a natural way.

\begin{pro}
Let $(M,\mu,\rho)$ be an $L$-structure and $\bar M$ its completion as a metric space.
Then there is a metric-measure $L$-structure on $\bar M$ such that $M\preceq\bar M$.
\end{pro}
\proof It is well-known that $\bar M$ carries a metric structure in the natural way.
What is new here is the measure part.
For each $L$-formula $\phi(\x)$, $\phi^M$ is uniformly continuous and hence it has a
unique continuous extension to $\bar M$. Let denote this extension by $\phi_{\bar M}$.
Let $\Ac$ be the $\sigma$-algebra of subsets of $\bar M$ generated by sets of the form
$\{x\in\bar M:\ \phi_{\bar M}(\a,x)\geqslant0\}$ where $\a\in M$.
For any $X\in\Ac$, define $\nu_0(X)=\mu(X\cap M)$. Note that this is a well-defined function
and indeed a probability measure on $\bar M$. Since open balls are included in this $\sigma$-algebra,
by \ref{Borel extension} and \ref{uniqueness}, $\nu_0$ extends to a unique Borel measure,
say $\nu$, on $\bar M$.
Now, $(\bar M,\nu,\rho_{\bar M})$ is a metric-measure $L$-structure.
We must show that it is an elementary extension of $M$.
For this purpose, it is sufficient to verify that for each $L$-formula
$\phi(\x)$, $\phi^{\bar M}=\phi_{\bar M}$.
We may do this by induction on the complexity of $\phi$. The main steps are quantifier cases.
For the integrations case, assume the claim holds for $\phi(\x,y)$.
First, for each $\a\in M$ we have
$$\phi^{\bar M}(\a)=\int_{\bar M}\phi^{\bar M}(\a,y)d\mu=
\int_{\bar M}\phi_{\bar M}(\a,y)d\mu=\int_{M}\phi^M(\a,y)d\nu.$$
For the last equality note that for each $r<s$, the set
$\{y\in\bar M: r<\phi_{\bar M}(\a,y)\leqslant s\}$ has the same measure as
its intersection with $M$. Since $M$ is dense in $\bar M$, by continuity,
they are also identical on $\bar M$. The case $\sup_y\phi(\x,y)$ is obvious. $\square$
\bigskip

Now we prove the downward theorem.
Let $(M,\mu)$ be a metric-measure structure in $L$ and $N$ be a metric substructure of $M$.
Let $\nu$ be the subspace measure on $N$. As stated before, $\nu$ is a Borel measure on $N$.
So, $(N,\nu)$ is a metric-measure structure in $L$.

\begin{pro}
Let $N$ be a  metric-measure structure in a countable language $L$ and $\kappa$ be
a cardinal such that $\kappa^{\aleph_0}=\kappa$.
Then for every $X\subseteq N$ with $|X|\leqslant\kappa$ there is an elementary
substructure $M\preceq N$ of cardinality $\kappa$ containing $X$.
\end{pro}
\proof Without loss of generality assume $X=M_0$ is a metric substructure of $N$ of
cardinality $\kappa$. We can easily define a countable chain
$$M_0\subseteq M_1\subseteq\cdots$$ of metric substructures of $N$ such that for each $n$,
$|M_n|=\kappa$ and

- for every $\phi(x)$ with parameters in $M_n$ and $\epsilon>0$,
there are $c\in M_{n+1}$ such that $\sup_x\phi^{M}(x)-\epsilon\leqslant\phi^{M}(c)$

- for each $X\subseteq N$ of positive measure in the $\sigma$-algebra generated
by formulas $\phi(x)$ with parameters in $M_n$,\ $X\cap M_{n+1}\neq\emptyset$.

To obtain the second clause we use the assumption that $\kappa^{\aleph_0}=\kappa$.
Let $M=\cup_n M_n$ and put the subspace measure $\mu_M$ on $M$ (as well as the metric substructure).
Then $(M,\mu_{M})$ is a metric-measure substructure of $N$.
Note that by definition every subset of $N$ in the $\sigma$-algebra generated by formulas with
parameters in $M$ has nonempty intersection with $M$. So, $M$ has full outer measure in $N$.
Now, to show that $M\preceq N$, we use induction.
Obviously, for each atomic $\phi(\x)$ and $\a\in M$ we have that $\phi^M(\a)=\phi^N(\a)$.
The connective cases and also the $\sup$ case are easy. Assume the claim holds for $\phi(\x,y)$.
Then by \ref{subspacemeasure} and induction hypothesis we have that
$$\int_N \phi^N(\a,y)dy=\int_M\phi^N(\a,y)dy=\int_M\phi^M(\a,y)dy=\Big(\int\phi\ dy\Big)^M(\a).\ \square$$

\subsection{Saturation and definability}
The notions of definable relation, definable set, type and saturation
are defined as in continuous logic. One can easily show that

\begin{pro}
Let $L$ be a countable language and $\mathcal D$ a countably incomplete
ultrafilter on $I$. Then $\prod_{\mathcal D} M_i$ is $\aleph_1$-saturated.
\end{pro}

Let $T$ be a complete theory and $(M,d,\mu)$ an $\aleph_0$-saturated model of $T$.
So, every type over the empty set is realized in $M$.
There are two kinds of topology on $S_n(T)$. The logic topology is generated
by complements of the sets of the form
$$[\phi\leqslant r]=\{p: p\vDash \phi\leqslant r\}.$$
The metric topology is generated by the metric
$$d(p,q)=\inf\{d(\a,\b):\ M\vDash p(\a),\ M\vDash q(\b)\}$$
$$d(\a,\b)=\max_i d(a_i,b_i).$$
The metric topology is stronger than the logic topology.
Moreover, $S_n(T)$ is compact with respect to the logic topology
and complete with respect to the metric topology.
Note that, for each $n$, the surjection $\a\mapsto tp(\a)$ from $M^n$ onto
$S_n(T)$ is continuous with respect to the metric topology on $S_n(T)$.
Hence, it induces a measure on $S_n(T)$ which by \ref{Borel extension} and
\ref{uniqueness} extends to a unique metric-measure on $S_n(T)$.
We denote this measure by $\mu_{S_n}$.

Given a formula $\phi(\x)$, for each type $p(\x)$ there is a unique $r$ such that
$\phi=r\in p$. We denote this $r$ by $\phi(p)$. The map $p\mapsto\phi(p)$
is denoted by $\bar\phi$. We recall a proposition from \cite{BBHU}.

\begin{pro} \label{Ben}
For any bounded function $f:S_n(T)\rightarrow\Rn$ the following are equivalent:

- $f$ is continuous for the logic topology

- there is a sequence $\phi_k(\x)$ of formulas such that $\bar\phi_k$ converges
to $f$ uniformly on $S_n(T)$

- $f$ is continuous for the logic topology and uniformly continuous for the $d$-metric.
\end{pro}

A natural question is what does happen if $f$ is assumed to be measurable.
The answer is easy in the $\aleph_0$-categorical case. Recall that if $T$ is
$\aleph_0$-categorical then logic topology and metric topology coincide on $S_n(T)$.

\begin{pro}
Assume $T$ is $\aleph_0$-categorical and $f:S_n(T)\rightarrow\Rn$
is a bounded function. Then $f$ is $\mu_{S_n}$-measurable if and
only if there is a sequence $\phi_k(\x)$ of formulas such that
$\bar\phi_k$ converges to $f$ pointwise on $S_n(T)$.
Moreover, in this case, $\int f dp=\lim_k\int\phi_k d\x$.
\end{pro}
\proof The `if' part is obvious. Let us prove the `only if' part.
It is known that if $X$ is a metric space and $\mu$ is a Borel measure on it,
then bounded continuous functions are dense in $L^1(X)$.
Since $f$ is bounded, it belongs to $L^1(S_n)$.
So, there is sequence $f_n$ of bounded continuous functions on $S_n(T)$
such that $\int|f-f_k|\searrow0$.
Note that this implies that $f_k\rightarrow f$ pointwise.
Now, by Proposition \ref{Ben}, we may assume without loss of generality
that $f_k$ is of the form $\bar\phi_k$ for some formula $\phi_k$.
So, $\bar\phi_k\rightarrow f$ pointwise. For the second part use
convergence theorems and the change of variable formula. $\square$

\section{Some applications} \label{applications}
In this section we give two applications of the compactness theorem.
Below, by ``$\phi(\x)=0$ for all $x$'' we mean the statement $\sup_{\x}|\phi(\x)|=0$.
The following is one of the various Riesz representation theorems (see \cite{parthasarathy}).

\begin{thm} \label{Riesz}
Let $M$ be a compact metric space and $I$ a positive linear functional
on $C(M)$ with $I(1)=1$. Then there exists a Borel measure $\mu$ on $M$
such that $I(f)=\int f d\mu$ for every $f\in C(M)$.
\end{thm}
\proof
By Dini's theorem, if a sequence $f_k\in C(M)$ decreases to zero pointwise, it converges
uniformly to zero and hence $I(f_k)$ tends to $0$. Indeed, $I$ is a Daniell-integral on $C(M)$.
Let $L$ be the language consisting of a constant symbol $c_a$ for each $a\in M$
and a unary relation symbol $R_f$ for each $f\in C(M)$.
Set $\flat_{R_f}=\sup_x |f(x)|$ and choose a modulus of uniform continuity
$\Delta_{R_f}$ with respect to which $f$ is uniformly continuous.
By the assumptions, we have that $I(|f|)\leqslant\sup_x|f(x)|$ for every $f$.
Let $T$ be the following $L$-theory where $a,b\in M$ and $f,g,h\in C(M)$:

1. $d(c_a,c_b)=d^M(a,b)$

2. $R_f(c_a)=f(a)$

3. $R_1(x)=1$\ for all $x$

4. $R_f(x)=rR_g(x)$\ \ for all $x$,\ \ where $f=rg$

5. $R_f(x)=R_g(x)+R_h(x)$\ \ for all $x$, \ where $f=g+h$

6. $R_f(x)=R_g(x)\wedge R_h(x)$\ \ for all $x$, \ where $f=g\wedge h$

7. $\int R_f\ dx=I(f)$\ \ for each $f\in C(M)$.
\bigskip

We first show that $T$ is finitely approximately satisfiable.
Let $T_0$ be a finite part of $T$ and $1=f_1,f_2,...,f_t$ be the list
of functions for which $R_{f_i}$'s appear in $T_0$.
There is no harm if we assume each $f_i$ is nonnegative (just add a
big positive real to them).
Fix $\epsilon>0$. We will define a measure on $M$ and
interpret $R_{f_i}$'s on $M$ such that the axioms of $T_0$ hold with error of
at most $\epsilon$. Let $J=[0,\alpha)$ contain the range of every $f_i$.
Let $[u_1,u_2),\ldots,[u_{s-1},u_s)$ be a partition of $J$
each one having length less than $\epsilon$.
Let $\Bc_0$ be the Boolean algebra generated by the sets $f_i^{-1}[u_j,u_{j+1})$
and $\mathcal P=\{P_1,...,P_\ell\}$ be its atoms. Each $P_k$ is then of
the form $\bigcap_{i=1}^{t}f_i^{-1}[u_{j_i},u_{j_i+1})$.
Let $\xi_i\geqslant0$ be a $\mathcal P$-simple function with $0\leqslant f_i-\xi_i<\epsilon$.
Interpret $R_{f_i}$ by $\xi_i$. Then all instances of the first six axioms
appearing in $T_0$ hold in $M$ by a good approximation.
For the last axiom, we need a finitary measure on $\Bc_0$.
By the lattice properties of $C(M)$, for each $P_k$ there is a sequence $h_{kn}\in C(M)$
increasing to $\chi(P_k)$ the support of each one being a subset of $P_k$.
Set $\lambda(P_k)=\lim_n I(h_{kn})$ and extend it to measure on $\Bc_0$
in the natural way.
If $\xi_i=\sum_{k=1}^\ell r_{ik}.\chi_{P_k}$ then let $\xi_{in}=\sum_{k=1}^{\ell}r_{ik}.h_{kn}$.
Clearly, $f_i-\xi_{in}\geqslant0$ for each $i$ and $n$.
Moreover, for each $x\in M$,\ $\xi_i(x)=\lim_n \xi_{in}(x)$.
One can easily check that $I((f_i-\xi_{in})\vee\epsilon)$ decreases to $\epsilon$
as $n$ tends to the infinity. Therefore
$$|I(f_i)-\int\xi_i\ d\lambda|=|I(f_i)-\sum_{k=1}^\ell r_{ik}.\lim_n I(h_{kn})|$$
$$=\lim_n I(f_i-\xi_{in})\leqslant\lim_n I((f_i-\xi_{in})\vee\epsilon)=\epsilon.$$
In particular, $|1-\lambda(M)|\leqslant\epsilon$ and by normalizing the
measure we get back to the probability case while retaining the required
approximations. This shows that $T$ is finitely approximately satisfiable.

Let $(N,\nu)$ be a model of $T$. We may suppose without loss of generality that
$M$ is a subset of $N$ and that each $f\in C(M)$ is the restriction to $M$ of $R_f^N$.
We first show that $M$ has full outer measure in $N$.
The outer measure of $M$ is equal to the infimum of the sums
$\sum_k\nu(A_k)$ where $M\subseteq\cup A_k$ and each $A_k$ is a Borel subset of $N$.
By Assumption \ref{assumption}, Proposition \ref{Bruckner} and Proposition \ref{uniqueness},
there is no harm if we assume the $A_k$'s are chosen from a smaller family of sets
containing a basis and closed under finite intersections.
In particular, since Borel and Baire $\sigma$-algebras coincide on $N$,
we may assume each $A_k$ is a finite intersection of sets of the form
$(R^N_{f})^{-1}(J)$ where $f\in C(M)$ and $J$ is an open interval.
Indeed, using the lattice properties, for each $k$ we may find
$f_k\in C(M)$ such that $A_k=(R^N_{f_k})^{-1}(0,\infty)$.
Then, by compactness of $M$, there exists $m$ such that $M\subseteq\cup_{k=1}^m A_k$.
Therefore, setting $f=\vee_{k=1}^m f_k$, we have that
$$M\subseteq\cup_{k=1}^m A_k=\{a\in N:\ R^N_f(a)>0\}.$$
Again, by compactness of $M$, there exists $r>0$ such that $f\geqslant r$ on $M$.
So, the axioms of $T$ imply that $R_f^N\geqslant r>0$ on $N$.
We conclude that $\sum_{k=1}^m\nu(A_k)\geqslant1$. So, $\nu^*(M)=1$.

Let $\mu$ be the subspace measure induced by $\nu$ on $M$ which is clearly a Borel measure.
Then by Proposition \ref{subspacemeasure} for each $f\in C(M)$ we have that
$$\int_M f d\mu=\int_N R_f^N d\nu=I(f).\ \square$$

The second application of compactness is the existence of invariant
measures on compact metric spaces.

\begin{thm}
Let $M$ be a nonempty compact metric space with isometry group $G$.
Then there exists a $G$-invariant Radon measure on $M$.
\end{thm}
\proof By \cite{Fremlin} 441 L, a measure $\mu$ on $M$ is $G$-invariant if and only if
for every $\alpha\in G$ and $f,g\in C(M)$ with $g=f\circ\alpha$ one has that $\int gd\mu=\int fd\mu$.
Let $L$ be the language consisting of a constant symbol $c_a$ for each $a\in M$ and
a unary relation symbol $R_f$ for each $f\in C(M)$.
Set $\flat_{R_f}=\sup_x |f(x)|$ and choose a modulus of uniform continuity
$\Delta_{R_f}$ with respect to which $f$ is uniformly continuous.
Let $T$ be the following $L$-theory where $a,b\in M$ and $f,g,h\in C(M)$:

1. $d(c_a,c_b)=d^M(a,b)$

2. $R_f(c_a)=f(a)$

3. $R_1(x)=1$\ for all $x$

4. $R_f(x)=rR_g(x)$\ \ for all $x$,\ \ where $f=rg$

5. $R_f(x)=R_g(x)+R_h(x)$\ \ for all $x$, \ where $f=g+h$

6. $R_f(x)=R_g(x)\wedge R_h(x)$\ \ for all $x$, \ where $f=g\wedge h$

7. $\int 1\ dx=1$

8. $\int R_g(x)\ dx=\int R_f(x)\ dx$\ \ if $\alpha\in G$ and $g=f\circ\alpha$.
\bigskip

Let $T_0$ be a finite part of $T$ and $\epsilon>0$ be fixed.
Assume $R_{f_1},...,R_{f_k}$ and $c_{a_1},...,c_{a_\ell}$ appear in $T_0$.
Let $\Delta$ be the minimum of moduli of continuities of $f_1,...,f_k$
and $2\delta>0$ be less than both $\Delta(\epsilon)$ and the minimum distance between $a_i$'s.

Let $B_{\delta}(b_1),...,B_{\delta}(b_n)$ be a set of balls of
radius $\delta$ covering $M$ and $B_{\delta}(c_1),...$, $B_{\delta}(c_m)$
be the least number of balls of radius $\delta$ containing $b_1,...,b_n$.
Then each $B_{\delta}(c_i)$ contains at least one of the $b_i$'s
and every $x\in M$ is in distance at most $2\delta$ of some $c_i$.
Let $X=\{c_1,...,c_m\}$ and put the uniform probability measure $\nu_0$ on it.
All instances of axioms 1-7 appeared in $T_0$ are satisfied in $C$ with an error
a multiple of $\epsilon$ which only depends on $T_0$ (does not depend on $\epsilon$).
Suppose the statement $\int R_f dx=\int R_g dx$ appears in $T_0$
where $g=f\circ\alpha$. Then
$$\int f|_X\ d\nu_0={1\over m}\sum_{i=1}^mf(c_i), \ \ \ \ \ \
\int g|_X\ d\nu_0={1\over m}\sum_{i=1}^mg(c_i).$$
By minimality of $|X|$, for each $c_i$
the ball $\alpha^{-1}(B_{\delta}(c_i))$ contains some $b_j$.
This $b_j$ is itself contained in some $B_{\delta}(c_{i'})$.
So, $d(c_{i'}, \alpha^{-1}(c_i))<2\delta$ and hence $d(\alpha(c_{i'}),c_i)<2\delta$.
More generally, for each distinct $c_{i_1},...,c_{i_k}$
there must exist distinct $c_{i'_1},...,c_{i'_k}$ such that
$$d(\alpha(c_{i'_t}), c_{i_t})<2\delta, \ \ \ \ \ t=1,...,k.$$
Indeed, if $c_{i'_1},...,c_{i'_{k'}}$ are the only points such that
every
$b_j\in B_{\delta}(\alpha^{-1}(c_{i_1}))\cup\cdots\cup B_{\delta}(\alpha^{-1}(c_{i_k}))$
is contained in $B_{\delta}(c_{i'_1})\cup\ldots\cup B_{\delta}(c_{i'_{k'}})$,
then $k'<k$ would contradict the minimality of $|X|$.
Using Hall's marriage theorem for bipartite graphs, one can arrange that
$c_i\mapsto c_{i'}$ be a permutation of $X$.
Then we have that
$$|g(c_{i'})-f(c_i)|=|f(\alpha(c_{i'}))-f(c_i)|<\epsilon$$
and hence $$\Big|\int f|_X\ d\nu_0-\int g|_X\ d\nu_0\Big|\leqslant\epsilon.$$
This shows that $T_0$ is finitely approximately satisfiable. Hence, $T$ is satisfiable.
Let $(N,\nu)$ be a model of $T$ which we may assume contains $M$ as a subset.
As in Theorem \ref{Riesz}, one easily shows that $M$ has full outer measure in $N$.
Let $\mu$ be the subspace measure on $M$.
Then, by Proposition \ref{subspacemeasure}, for each $f,g\in C(M)$
and $\alpha\in G$, if $g=f\circ\alpha$ one has that
$$\int_M g\ d\mu=\int R_g^N d\nu=\int R_f^N d\nu=\int_M f\ d\mu.\ \square$$

\end{document}